\documentclass[a4paper,12pt]{article}
\usepackage{amssymb,amsmath ,bm}
\def\beq{\begin{equation}}
\def\eeq{\end{equation}}
\def\bea{\begin{eqnarray}}
\def\eea{\end{eqnarray}}
\def\nn{\nonumber}
\def\U{{\cal U}}
\def\A{{\cal A}}

\def\hf{\frac{1}{2}}
\def\bra#1{\left\{ #1 \right\}_{q^{-1}}}
\def\ku#1{\left\{ #1 \right\}_q}
\def\qn#1{[#1]_q}
\def\ee#1{e^{\ell}_{#1}(\lambda)}
\def\h#1{\hat{#1}}
\def\sq#1{\big\lgroup #1 \big\rgroup_q}
\def\sqi#1{\big\lgroup #1 \big\rgroup_{q^{-1}}}
\def\uT#1{{\cal T}_{#1}}
\def\Tm#1#2{T^{#1}_{#2}}
\def\e#1#2{e^{#1}_{#2}}
\def\CGC#1#2{C^{#1}_{#2}\,}
\def\Pz#1#2{P^{#1}_{#2}(\zeta)}

\begin{document}
%
%
%
%
\thispagestyle{empty}

\vspace*{3cm}

\begin{center}
{\LARGE\sf
Basic Hypergeometric Functions and Covariant Spaces 
for Even Dimensional Representations of \bm{$U_q[osp(1/2)]$}
}

\bigskip\bigskip
N. Aizawa

\bigskip
\textit{
Department of Mathematics and Information Sciences, \\
Graduate School of Science,\\
Osaka Prefecture University, \\
Nakamozu Campus, Sakai, Osaka 599-8531, Japan}\\

\bigskip
R. Chakrabarti, S.S. Naina Mohammed

\bigskip
\textit{
Department of Theoretical Physics, \\
University of Madras, \\
Guindy Campus, Chennai 600 025, India
}

\bigskip
J. Segar

\bigskip
\textit{
Department of Physics, \\
Ramakrishna Mission Vivekananda College, \\
Mylapore, Chennai 600 004, India
}

\end{center}

\vfill
\begin{abstract}
Representations of the quantum superalgebra $ U_q[osp(1/2)] $ and their 
relations to the basic hypergeometric functions are investigated. We 
first establish  Clebsch-Gordan decomposition for the superalgebra 
$ U_q[osp(1/2)] $ in which the representations having no classical 
counterparts are incorporated. Formulae for these Clebsch-Gordan 
coefficients are derived, and it is observed that they may be expressed 
in terms of the $Q$-Hahn polynomials. We next investigate 
representations of the quantum supergroup $ OSp_q(1/2) $ which are not 
well-defined in the classical limit. Employing the universal 
$ {\cal T}$-matrix, the representation matrices are obtained explicitly, 
and found to be related to the little $Q$-Jacobi polynomials. 
Characteristically, the relation $ Q = -q $ is satisfied in all cases. 
Using the Clebsch-Gordan coefficients derived here, we construct new 
noncommutative spaces that are covariant under the coaction of the 
even dimensional representations of the quantum supergroup $OSp_q(1/2)$.
\end{abstract}


\newpage

\setcounter{page}{1}
%
%
%
\section{Introduction}
\label{Intro}

  Soon after the introduction of the quantum groups, their relation to 
the basic hypergeometric functions via the representation theory was 
revealed by many authors \cite{KS}. In particular, it was observed 
\cite{VS,Koo,MMNNU} that the representation matrices of the quantum 
group $ SU_q(2) $ can be expressed in terms of the little $q$-Jacobi 
polynomials. Kirillov and Reshetikhin demonstrated \cite{KiRe} that  
the Clebsch-Gordan coefficients of the quantum algebra $ U_q[sl(2)] $ 
relate to the $q$-Hahn and the dual $q$-Hahn polynomials. These works 
furnished a new algebraic framework to the theory of basic 
hypergeometric functions. Since then extensive studies on interrelations 
between the quantum group representations and the basic hypergeometric 
functions had taken place. We here mention some key examples of these 
developments. The matrix elements of  the quantum group $ SU_q(1,1) $ were 
found to be related \cite{MMNMSU} to the polynomials obtained from 
$ {}_2\phi_1$. The realizations of the quantum algebra $ U_q[su(1,1)] $ 
and the generating functions of Al-Salam-Chihara polynomials were found  
\cite{JJ97} to be linked. The kinship between the group theoretical 
treatment of the $q$-oscillator algebra and the $q$-Laguerre as well as 
the $q$-Hermite polynomials was observed \cite{FV91}. The 
connection between the metaplectic representation of $ U_q[su(1,1)] $ 
and the $q$-Gegenbauer polynomials was noted \cite{FV92}.

  On the other hand, study of relations between the quantum supergroups 
and the basic hypergeometric functions started very recently. In 
\cite{ZZ} the homogeneous superspaces for the general linear supergroup 
and the spherical functions on them were investigated. Zou studied 
\cite{Zou06} the spherical functions on the symmetric spaces arising 
from the quantum superalgebra $U_q[osp(1/2)]$. This author also observed 
\cite{Zou03} the relationship between the transformation groups of the 
quantum super 2-sphere and the little $Q$-Jacobi polynomials. Considering a 
$2 \times 2$ quantum supermatrix and identifying its dual algebra with 
the quantum superalgebra $U_q[osp(1/2)]$, the finite dimensional 
representations of the quantum supergroup $ OSp_q(1/2) $ were found 
\cite{Zou03} to be related to the little $Q$-Jacobi polynomials with the 
assignment $ Q = -q$ . Instead of $Q,$ the parameter $ t = i\sqrt{q} $ 
was used in \cite{Zou03}. 
Adopting an alternate procedure by explicitly 
evaluating the universal ${\cal T}$-matrix that capped the Hopf dual 
structure, and using the representations of the $ U_q[osp(1/2)] $ 
algebra, the present authors obtained \cite{ACNS} the same result  
independently. The results in \cite{ACNS} are, however, partial in the 
sense that only odd dimensional representations of the algebra 
$U_q[osp(1/2)]$ are taken into account. One of the purposes of the 
present study is to incorporate even dimensional representations of the 
supergroup $OSp_{q}(1/2)$ in the framework of \cite{ACNS}. Continuing 
our study of the finite dimensional representations of the universal 
${\cal T}$-matrix, we observe that the even dimensional representations 
of the quantum supergroup $OSp_{q}(1/2)$ may too be expressed via the  
little $Q$-Jacobi polynomials with $ Q = -q $. Furthermore, we also 
study irreducible decomposition of the tensor product of both even and 
odd dimensional representations of the $U_q[osp(1/2)]$ algebra. 
Evaluating the Clebsch-Gordan coupling of two even dimensional 
representations as well as that of an even and an odd dimensional 
representations, we observe that the decomposition is multiplicity free. 
Proceeding further we notice that the Clebsch-Gordan coefficients for the 
decompositions are related to the $Q$-Hahn polynomials with $Q = -q$.

 Emergence of the $ Q = -q $ polynomials, in contrast to the $ Q = q $ 
polynomials being present for the aforementioned quantum groups, appears 
to be a generic property of the quantum superalgebra $ U_q[osp(1/2)].$   
For odd dimensional representations, this property may be interpreted as 
a reflection of the isomorphism of $ U_q[osp(1/2n)] $ and 
$ U_{-q}[so(2n+1)] $ which holds on the representation space 
\cite{Zhang}. The even dimensional representations for which 
the said isomorphism is not known, however, are still characterized by   
polynomials with $ Q = -q $. Pointing towards a generalized 
feature of the quantum {\it supergroups} the present work puts forward 
new entries in the list of relations between supergroups and basic 
hypergeometric functions.

 Explicit evaluation of the Clebsch-Gordan coefficients allows us 
to explore {\it new noncommutative spaces covariant under the action of 
even dimensional representations of} $ OSp_q(1/2). $ Employing the method 
developed in \cite{AC}, we, for instance, introduce the defining 
relations of the covariant noncommutative space of dimension four. 
Our construction may be generalized to describe similar covariant 
noncommutative spaces of higher dimensions. Especially for the root of 
unity values of $q$ the representation of these spaces may be of 
interest in some physical problem.

Our focus on $OSp_{q}(1/2)$ is explained by physical and mathematical 
reasons: Physically, fully developed representation theory of the said 
supergroup may provide better insight to the solvable vertex type models 
\cite{Sal} endowed with the quantum $U_{q}[osp(1/2)]$ symmetry, Gaudin 
models \cite{KM} and two dimensional field theories \cite{SK}. 
Mathematical motivation lies in the fact that $ osp(1/2) $ is the 
simplest superalgebra and a basic building block for other superalgebras.

We plan the paper as follows. In the next section, definitions and 
representations of $ U_q[osp(1/2)]$ to be used in the subsequent 
sections, are listed. We prove that the even dimensional representations 
are of grade star type. Tensor product of two irreducible 
representations is considered in \S \ref{CGC}. We show that the tensor 
product is decomposed into direct sum of irreducible representations 
without multiplicity. Formulae of the Clebsch-Gordan coefficients 
are derived and it is shown that they are expressed in terms of the 
$Q$-Hahn polynomials with $ Q = -q.$ In \S \ref{RepA}, even dimensional 
representations of the quantum supergroup $ OSp_q(1/2) $ are constructed 
using the universal $ {\cal T}$-matrix. Relations of the matrix elements 
and the little $Q$-Jacobi polynomials with $ Q = -q $ is established. 
We discuss covariant noncommutative superspaces of dimensions two and 
four in \S \ref{CovSpa}. Our concluding remarks are given 
in \S \ref{ConclR}.

%
%
%

\setcounter{equation}{0}
\section{\bm{$U_q[osp(1/2)]$} and its representations}
\label{UqRep}

The quantum superalgebra $ \U \equiv U_q[osp(1/2)] $ has been introduced 
in \cite{Kul}. Finite dimensional representations of $ \U, $ which are 
$q$-analogue of the representations of the classical superalgebra 
$ osp(1/2), $ have been investigated in \cite{KR,Sal}. Classification of 
finite dimensional integrable representations of more general quantum 
superalgebra $ U_q[osp(1/2n)] $ has been made by Zou \cite{Zou}, and 
therein it has been observed that $U_{q}[osp(1/2)]$ admits representations 
which are not deformation of the ones for $ osp(1/2)$. Thus we have two 
types of finite dimensional representations for $ \U: $ one of them has 
classical counterparts, and the other does not. For the purpose of 
fixing our notations and conventions we here list the relations that will 
be used subsequently.

  The algebra $\U$ is generated by three elements $H$ (parity even) and 
$ V_{\pm} $ (parity odd) subject to the relations
\beq
   [H, V_{\pm}] = \pm\hf V_{\pm}, \qquad 
   \{ V_+, V_- \} = -\frac{q^{2H} - q^{-2H}}{q-q^{-1}} \equiv -\qn{2H}.
   \label{ospS2}
\eeq
The deformation parameter $q$ is assumed to be generic throughout this article. 
The Hopf algebra structures defined via the coproduct ($\Delta$), 
the counit ($\epsilon$) and the antipode ($S$) maps read as follows:
\bea
  & & \Delta(H) = H \otimes 1 + 1 \otimes H, \qquad
      \Delta(V_{\pm}) = V_{\pm} \otimes q^{-H} + q^H \otimes V_{\pm},
      \label{Dosp2} \\
  & & \epsilon(H) = \epsilon(V_{\pm}) = 0, \label{epS2} \\
  & & S(H) = -H, \qquad S(V_{\pm}) = -q^{\mp 1/2} V_{\pm}. \label{SospS2}
\eea

  The finite dimensional irreducible representations of $ \U $ are 
specified by the highest weight $ \ell $ which takes any non-negative 
integral or half-integral value. We denote the irreducible 
representation space  of highest weight $\ell$ by $ V^{(\ell)}$. 
According to \cite{Zou}, 
a representation is referred to be integrable if $ V^{(\ell)} $ 
is a direct sum of its weight spaces, and if $ V_{\pm}$ act as locally 
nilpotent operators on $ V^{(\ell)}.$ The results in \cite{Zou} for the 
case of $ \U $ can be stated quite simply by introducing an element 
$ K = q^{2H} $ as follows: 
Let $ \bm{v} \in V^{(\ell)} $ be a highest weight vector ($ V_+ \bm{v} = 0$). 
The highest weight representation constructed on $ \bm{v} $ is 
integrable if and only if
\beq
    K \bm{v} = 
    \left\{
        \begin{array}{lcl}
            \pm q^{\ell} & \quad & \mbox{if $ \ell $ is an integer}, \\
            \pm i q^{\ell} &     & \mbox{if $ \ell $ is a half-integer},
        \end{array}
    \right.
    \label{HWMod}
\eeq
and the integrable representations are completely reducible. 
The representation in $ V^{(\ell)} $ has dimension $ 2 \ell + 1 $ 
so that $ V^{(\ell)} $ is odd (resp. even) dimensional if $ \ell $ is 
an integer (resp. half-integer). 
It is known that the classical superalgebra $ osp(1/2) $ does not have 
even dimensional irreducible representations.

  We denote a basis set of $ V^{(\ell)} $ as 
$ \{\; e^{\ell}_m(\lambda) \ | \ m = \ell, \ell-1,\cdots, -\ell \;\}, $
where the index $ \lambda = 0, 1 $ specifies the parity of the highest weight 
vector $ e^{\ell}_{\ell}(\lambda). $ The parity of the vector 
$ e^{\ell}_m(\lambda) $ equals $ \ell - m + \lambda,$ 
as it is obtained by the action of $ V_-^{\ell-m} $ on 
$ e^{\ell}_{\ell}(\lambda).$  For the superalgebras the norm of the 
representation basis need not be chosen positive definite.  
In this work, however, we assume the positive definiteness of 
the basis elements: 
\beq
 (e^{\ell}_m(\lambda), e^{\ell'}_{m'}(\lambda)) = \delta_{\ell \ell'} 
 \delta_{mm'}.
 \label{Normalization}
\eeq
With these settings, the irreducible representation of $ \U$ on $ V^{(\ell)} $ 
is given as follows: 
For $ \ell $ integer, we take the following form which is a variant of 
the convention used in \cite{KR,Sal}:
\bea
  & & H \e{\ell}{m}(\lambda) = \frac{m}{2}\, \e{\ell}{m}(\lambda), \nn \\
  & & V_+ \e{\ell}{m}(\lambda) = \left(\frac{1}{\ku{2}}\ku{\ell-m}\ku{\ell+m+1}\right)^{1/2} 
      \e{\ell}{m+1}(\lambda),
      \label{RepU} \\
  & & V_- \e{\ell}{m}(\lambda) 
     = (-1)^{\ell-m-1} \left(\frac{1}{\ku{2}} \ku{\ell+m}\ku{\ell-m+1}\right)^{1/2}
      \e{\ell}{m-1}(\lambda), \nn
\eea
where
\beq
  \ku{m} = \frac{q^{-m/2}-(-1)^m q^{m/2}}{q^{-1/2}+q^{1/2}}.  \label{Ksymb}
\eeq
The representation space $ V^{(\ell)} $ is odd dimensional. 
It is known that (\ref{RepU}) is a grade star representation \cite{SNR} 
if $ q \in {\mathbb R} $ \cite{Sal}. The grade adjoint operation is given by
\beq
  H^* = H, \qquad 
  V_{\pm}^* = \pm (-1)^{\epsilon} V_{\mp}, \label{Gad_odd}
\eeq
where $ \epsilon = \lambda + 1 $ (mod 2). 
The grade adjoint operation is assumed to be an algebra anti-isomorphism and a 
coalgebra isomorphism.

  For a half-integer $ \ell $, we chose a representation parallel to 
(\ref{RepU}) but different from the one in \cite{Zou}:
\bea
  & & H \ee{m} = \frac{1}{2}(m \pm \eta)\, 
  \ee{m}, 
  \nn \\
  & & V_+ \ee{m} = \pm \left(
    \frac{1}{ \ku{2} } \ku{\ell-m} \ku{\ell+m+1}
  \right)^{1/2} \ee{m+1}, 
  \label{RepU_even} \\
  & & V_- \ee{m} = (-1)^{\ell-m} i \left(
    \frac{1}{ \ku{2} } \ku{\ell+m} \ku{\ell-m+1}
  \right)^{1/2} \ee{m-1},
  \nn
\eea  
where
\[
  \eta = \frac{\pi i}{2 \ln q}.
\]
The factor $i$  appearing in the action of $K$ (\ref{HWMod}) is converted 
into the constant $ \eta $ for the action of $H.$ 
We keep two different phase conventions for later convenience. 
Representation spaces corresponding to each phase choice are 
denoted by $ V^{(\ell)}_{\pm}. $ 
The space $ V^{(\ell)}_{\pm} $ is even dimensional. 
We come to state our first result. If $ q \in {\mathbb R},$ then 
the (\ref{RepU_even}) is a grade star representation under the 
grade adjoint operation
\beq
    H^{*} = H \mp \frac{\pi i}{2 \ln q}, \qquad
     V_+^{*} = \pm i(-1)^{\epsilon} \, V_-, \qquad
     V_-^{*} = \mp i(-1)^{\epsilon} \, V_+,
     \label{Gad_even}
\eeq
where $ \epsilon = \lambda + 1 $ (mod 2). 

  Proof of this statement is rather straightforward. Recall the 
definition of grade star representation. Let $ \rho $ be a 
representation of a quantum superalgebra $U.$  Denoting a grade star 
operation by $*$ that is defined in $U$, 
the representation $ \rho $ is referred to be of grade star 
type if $ \rho(X)\,\forall X \in U $ satisfies 
\beq
  \rho(X^*) = \rho(X)^*, 
  \label{GstarRep}
\eeq
$ \rho(X)^* $ is the superhermitian 
conjugate defined by 
\[
  \rho(X)^*_{ij} = (-1)^{(\h{i} + \h{X})(\h{i}+\h{j})} \overline{\rho(X)}_{ji},
\]
where $ \h{a} $ denotes the parity of the object $a$.  It is easily 
verified that (\ref{RepU_even}) and (\ref{Gad_even}) satisfy 
(\ref{GstarRep}). 

%
%
%
\setcounter{equation}{0}
\section{Clebsch-Gordan decomposition and \bm{$Q$}-Hahn polynomials}
\label{CGC}

  In this section, we consider tensor product of two irreducible 
representations of $ \U$. As noted earlier, the algebra $\U$ has two 
types of grade star representations, namely odd and even dimensional 
ones. The former maintains one-to-one correspondence to representation 
of $ osp(1/2)$ of same dimensionality, whereas the latter has no classical 
analogue. We may consider three cases of tensor product, namely, the 
product of two odd dimensional representations, two even ones, and an odd 
and an even ones. We observe that for all the above three cases the  
tensor product of two irreducible representations is, in general, 
reducible, and may be decomposed into a direct sum of irreducible ones 
without multiplicity:
\beq
  V^{(\ell_1)} \otimes V^{(\ell_2)} = 
  V^{(\ell_1+\ell_2)} \oplus V^{(\ell_1+\ell_2 -1)} \oplus \cdots 
  \oplus V^{(|\ell_1-\ell_2|)}.
  \label{CGdecomp}
\eeq
The decomposition of the tensored vector space in the irreducible basis 
is provided by the Clebsch-Gordan coefficients (CGC):
\beq
    \e{\ell}{m}(\ell_1,\ell_2,\Lambda) = 
    \sum_{m_1,m_2\atop m_{1}+m_{2} = m} 
    \CGC{\ell_1\ \ell_2\ \;\ell}{m_1\,m_2\,m}\,\,
    \e{\ell_1}{m_1}(\lambda) \otimes \e{\ell_2}{m_2}(\lambda),
    \label{CGCdef}
\eeq
where $\Lambda = \ell_1+\ell_2-\ell \ ({\rm mod}\ 2)$ signifies the 
parity of the highest weight vector 
$ \e{\ell}{\ell}(\ell_1,\ell_2,\Lambda).$ The decomposition 
(\ref{CGdecomp}) is established below by explicit construction of the 
CGC. We also demonstrate that the tensor product representation is not 
of grade star type. Another pertinent problem of interest is the 
interrelation of the CGC and the basic hypergeometric functions. We
prove below that the CGC have polynomial structure corresponding to the 
$Q$-Hahn polynomials. We treat three cases separately. 

%
%
\subsection{Two odd dimensional representations}
\label{OO}

The Clebsch-Gordon decomposition for $\U$ in this case has been  
extensively studied in \cite{MM}. We here discuss a relation of the CGC 
and the $Q$-Hahn polynomials. 

  Following \cite{GR} the basic hypergeometric function 
$ {}_{r+1}\phi_r $ is defined as 
\beq
  {}_{r+1}\phi_r\left[ \left. 
    \begin{array}{c}
     a_1, a_2, \dots, a_{r+1} \\ b_1, b_2, \dots, b_r
    \end{array} 
    \right| 
    Q; z
  \right] 
  = 
  \sum_{k=0}^{\infty} 
   \frac{ (a_1;Q)_k (a_2;Q)_k \cdots (a_{r+1};Q)_k }
        { (b_1;Q)_k (b_2;Q)_k \cdots (b_r;Q)_k }
   \,\,\frac{z^k}{ (Q;Q)_k },
   \label{BHF}
\eeq
where the shifted factorial reads 
\beq
  (x;Q)_k = \left\{
     \begin{array}{ccl}
        1  & & \hbox {for} \,k = 0, \\
        \displaystyle \prod_{j=0}^{k-1} (1-xQ^j)  & & 
        \hbox {for} \,k \neq 0.
     \end{array}
  \right.
  \label{shift_fac}
\eeq
The $Q$-Hahn polynomials are defined \cite{GR} via $ {}_3\phi_2 $ in 
a standard way: 
\beq
  {\cal Q}_M(x;a,b,N;Q) = {}_3\phi_2\left[ \left.
    \begin{array}{c}
      Q^{- M}, abQ^{M+1}, Q^{-x} \\
      aQ, Q^{-N}
    \end{array}
    \right| 
    Q; Q
    \right],\quad M \le N.
    \label{QHahn}
\eeq
Setting $ a = Q^{\alpha}, b= Q^{\beta}, $ we obtain the following form 
of the $Q$-Hahn polynomials:
\beq
 {\cal Q}_M(x;\alpha,\beta,N;Q) = 
 \sum_{k} \frac{(Q^{-M};Q)_k (Q^{\alpha+\beta+M+1};Q)_k (Q^{-x};Q)_k}
               { (Q^{\alpha+1};Q)_k (Q^{-N};Q)_k}
 \,\,\frac{Q^k}{(Q;Q)_k}.
 \label{QHahn2}
\eeq

  Explicit formulae of the CGC are found in \cite{MM} and \cite{AC}. 
Up to a multiplicative factor that is irrelevant to the present 
discussion, the CGC reads
\bea
  \CGC{\ell_1\ \ell_2\ \;\ell}{m_1\,m_2\,m} 
  &=& N_1(\ell_1,\ell_2,\ell,m ; q) 
   \sum_{m_1+m_2=m} (-1)^{(\ell_1-m_1) \lambda + \hf (\ell_1-m_1) (\ell_1-m_1+1)}   
  \nn \\
  &\times&
      q^{-\hf m_1(m+1)}
    \left(  
      \frac{ \ku{\ell_1-m_1}! \ku{\ell_2-m_2}! }{ \ku{\ell_1+m_1}! \ku{\ell_2+m_2}! }
    \right)^{1/2}
  \nn \\
  &\times& \sum_k 
   (-1)^{k(\ell_1+\ell_2-m) + \hf k(k-1)} 
   q^{\hf k (\ell+m+1)}
  \nn \\
  &\times& 
   \frac{ \ku{\ell_1+\ell-m_2-k}! \ku{\ell_2+m_2+k}! }
        { \ku{\ell_1-\ell+m_2+k}! \ku{\ell_2-m_2-k}! \ku{\ell-m-k}! \ku{k}! },
  \label{CGC_OO}
\eea
where the index $k$ runs over all non-negative integers maintaining the 
argument of $ \ku{x} $ non-negative. The part of summation over $k$ can 
be regarded as a polynomial appearing in the CGC. To relate this to the 
$Q$-Hahn polynomial with $ Q = -q$, we recast the factorials in 
(\ref{CGC_OO}) as shifted factorials:
\bea
  & & \frac{ \ku{A}! }{ \ku{A+k}! } = q^{\frac{1}{4} k (2A+k-1)}
      \frac{(1+q)^k}{ (Q^{A+1};Q)_k}, 
  \label{Q-fact1}\\
  & & \frac{ \ku{A}! }{ \ku{A-k}! } = 
      (-1)^{\hf k(2A+3-k)} q^{\frac{1}{4} k(2A+3-k)} 
      \frac{(Q^{-A};Q)_k}{(1+q)^k},
  \label{Q-fact2}
\eea
where $ Q = -q$. Apart from a multiplicative constant, the summation 
over $k$ in (\ref{CGC_OO}) is now related to the polynomial form
(\ref{QHahn2}) with the assignment $Q = - q$:
\bea
  \sum_k \; ( \cdots ) &=& 
    \frac{ \ku{\ell_1+\ell-m_2}! \ku{\ell_2+m_2}! }
         { \ku{\ell_1-\ell+m_2}! \ku{\ell_2-m_2}! \ku{\ell-m}! } 
  \nn \\
  & & \times 
   \sum_{k} \frac{(Q^{-\ell_2+m_2};Q)_k (Q^{\ell_2+m_2+1};Q)_k (Q^{-\ell+m};Q)_k }
               { (Q^{\ell_1-\ell+m_2+1};Q)_k (Q^{-\ell_1-\ell+m_2};Q)_k }
   \frac{Q^k}{(Q;Q)_k}. 
   \label{sumk_OO}
\eea
The parameters of the $Q$-Hahn 
polynomial read:
\beq
  \begin{array}{lcl}
     \alpha = -\ell+\ell_1+m_2, & \ & \beta = \ell-\ell_1+m_2 \\
     N = \ell + \ell_1 -m_2, & & x = \ell-m, \qquad M = \ell_2-m_2.
  \end{array}
  \label{Hahn_label}
\eeq

%
%
\subsection{Two even dimensional representations}
\label{EE}

 It is expected that in this case the decomposed irreducible spaces 
have odd dimensions. We consider representation in the tensored space 
$ V^{(\ell_1)}_+ \otimes V^{(\ell_2)}_- $ so that the constant $ \eta $ 
appearing in the eigenvalues of $ H$ is eliminated. The CGC for this 
case may be computed in the standard way as outlined below. 

We start by determining the highest weight states in the direct product 
space $V^{(\ell_1)}_+ \otimes V^{(\ell_2)}_-$. 
A highest weight state is a linear 
combination of the basis of $ V^{(\ell_1)}_+ \otimes V^{(\ell_2)}_-: $
\beq
  e^{\ell}_{\ell}(\ell_1,\ell_2, \Lambda) = \sum_{m_1,m_2} C_{m_1,m_2} 
  e^{\ell_1}_{m_1}(\lambda) \otimes e^{\ell_2}_{m_2}(\lambda). 
  \label{ell-Lambda}
\eeq
The defining equations for the highest weight state read
\beq
  \Delta(H)\, e^{\ell}_{\ell}(\ell_1,\ell_2,\Lambda) 
  = \frac{\ell}{2} \, e^{\ell}_{\ell}(\ell_1,\ell_2,\Lambda), 
  \qquad
  \Delta(V_+)\, e^{\ell}_{\ell}(\ell_1,\ell_2,\Lambda) = 0.
  \label{highest-L}
\eeq
The first equality in (\ref{highest-L}) puts a constraint on the 
summation in (\ref{ell-Lambda}). Implying the identity
\[
  \Delta(H)\, e^{\ell}_{\ell}(\ell_1,\ell_2, \Lambda) = 
  \sum_{m_1,m_2} \hf(m_1+m_2) C_{m_1,m_2} 
  e^{\ell_1}_{m_1}(\lambda) \otimes e^{\ell_2}_{m_2}(\lambda)
  = 
  \frac{\ell}{2} e^{\ell}_{\ell}(\ell_1,\ell_2,\Lambda),
\]
it necessitates that the summation must obey the constraint 
$ m_1 + m_2 = \ell.$ As both $m_1$ and $m_2$ are half-integers, 
$ \ell $ takes integral values. The second equality in 
(\ref{highest-L}) produces a recurrence relation:
\bea
  & & q^{-\hf m_2} \sqrt{ \ku{\ell_1-m_1} \ku{\ell_1+m_1+1} }
      C_{m_1,m_2}
    \nn \\
  & &  \qquad
      = (-1)^{\ell_1-m_1+ \lambda+1}  q^{\hf(m_1+1)}
      \sqrt{ \ku{\ell_2+m_2} \ku{\ell_2-m_2+1} } C_{m_1+1,m_2-1}.
    \label{rec-EE}
\eea
The recurrence relation may be solved explicitly: 
\bea
  C_{m_1,m_2} &=& (-1)^{(\ell_1-m_1) \lambda + \hf (\ell_1-m_1) 
  (\ell_1-m_1-1)}
  q^{\hf (\ell+1)(\ell_1-m_1)}
  \nn \\
  &\times&
  \left(
    \frac{ \ku{\ell_1+\ell_2-\ell}! \ku{\ell_1+m_1}! \ku{\ell_2+m_2}! }
         { \ku{2\ell_1}! \ku{-\ell_1+\ell_2+\ell}! \ku{\ell_1-m_1}! 
           \ku{\ell_2-m_2}!}
  \right)^{1/2} C_{\ell_1,\ell-\ell_1}.
  \label{C-sol}
\eea
The highest weight vector has been determined uniquely up to an overall 
factor $ C_{\ell,\ell-\ell_1}$ that may be obtained by using the 
normalization. As it is unnecessary for our purpose, we leave the 
constant undetermined.  

  Other states in $ V^{(\ell_1)}_{+} \otimes V^{(\ell_2)}_{-} $ 
are obtained by repeated application of $ \Delta(V_-) $ on 
the highest weight state:
\beq
  \Delta(V_-)^{\ell-m} e^{\ell}_{\ell}(\ell_1,\ell_2,\Lambda) = 
  \sum_{m_1+m_2=\ell} C_{m_1,m_2} \Delta(V_-)^{\ell-m} 
  e^{\ell_1}_{m_1}(\lambda) \otimes e^{\ell_2}_{m_2}(\lambda).
  \label{state-EE}
\eeq
The state (\ref{state-EE}) may be easily recognized as an eigenstate of 
the operator $ \Delta(H) $ with the eigenvalue $ m/2.$ To express the 
state (\ref{state-EE}) as a linear combination of 
$ \e{\ell_1}{m_1}(\lambda) \otimes \e{\ell_2}{m_2}(\lambda), $ we expand 
$  \Delta(V_-)^{\ell-m} $ by using the binomial theorem for anti-commuting 
objects. For $q$-anti-commuting operators subject to  
$ q AB + BA = 0$, the expansion reads  
\beq
  (A+B)^n = \sum_{k=0}^n q^{\hf k(n-k)} 
  \frac{ \ku{n}! }{ \ku{k}! \ku{n-k}! } A^{k} B^{n-k}.
  \label{binomial}
\eeq
Setting $ A = q^H \otimes V_-, \ B = V_- \otimes q^{-H}$, we apply the 
expansion (\ref{binomial}) in (\ref{state-EE}). Following a redefinition 
of the summation variables, we obtain the expression of CGC given below:
\bea
  \CGC{\ell_1\ \ell_2\ \;\ell}{m_1\,m_2\,m} 
  &=& N_2(\ell_1,\ell_2,\ell,m ; q) 
   \sum_{m_1+m_2=m} (-1)^{(\ell_1-m_1) \lambda + \hf (\ell_1-m_1) (\ell_1-m_1-1)}   
  \nn \\
  &\times&
      q^{-\hf m_1(m+1)}
    \left(  
      \frac{ \ku{\ell_1-m_1}! \ku{\ell_2-m_2}! }{ \ku{\ell_1+m_1}! \ku{\ell_2+m_2}! }
    \right)^{1/2}
  \nn \\
  &\times& \sum_k 
   (-1)^{k(\ell_1+\ell_2-m) + \hf k(k+1)} 
   q^{\hf k (\ell+m+1)}
  \nn \\
  &\times& 
   \frac{ \ku{\ell_1+\ell-m_2-k}! \ku{\ell_2+m_2+k}! }
        { \ku{\ell_1-\ell+m_2+k}! \ku{\ell_2-m_2-k}! \ku{\ell-m-k}! \ku{k}! }.
  \label{CGC_EE}
\eea
This CGC is almost same as (\ref{CGC_OO}) discussed in the previous 
subsection, except for a sign factor that originates from the difference 
in phases between the odd and the even dimensional representations, 
given in (\ref{RepU}) and (\ref{RepU_even}), respectively. Moreover, 
the sign difference in the factors comprising of the sum over the 
index $k$ disappears when the expression in (\ref{CGC_EE}) is 
recast in terms of the shifted factorials. This leads to identical   
sums in (\ref{sumk_OO}) and (\ref{CGC_EE}) on the index $k$. We may,
therefore, immediately conclude that the CGC in (\ref{CGC_EE}) are 
related to the $Q$-Hahn polynomials with $Q=-q$, and the values of the 
parameters are given by (\ref{Hahn_label}). 

    The above construction of the eigenstates of $ \Delta(H) $ is just 
the standard procedure of highest weight construction leading to a 
multiplet of $ 2 \ell + 1 $ states from 
$ e^{\ell}_{\ell}(\ell_1,\ell_2,\Lambda) $ by repeated actions of 
$ \Delta(V_-).$ 
The $ 2\ell + 1 $ states are linearly independent, since they are the eigenvectors 
of $ \Delta(H) $ with different eigenvalues. Therefore, they form a basis of 
an invariant subspace in $ V^{(\ell_1)}_+ \otimes V^{(\ell_2)}_-. $   
It is an easy task to verify that eigenstates (\ref{state-EE}) belonging to different 
values of $ \ell $ are linearly independent. 
It also follows that, as  $ \ell \geq 0, $ its 
possible values are
$ \ell_1+\ell_2, \ell_1+\ell_2-1, \cdots, |\ell_1- \ell_2|.$  
The total number of the eigenstates of $ \Delta(H) $, of course,    
coincides with the dimension of $ V^{(\ell_1)}_{+} \otimes 
V^{(\ell_2)}_{-}$:
\[
  \sum_{\ell=|\ell_1-\ell_2|}^{\ell_1+\ell_2} (2 \ell+1)
  = (2\ell_1 + 1)(2 \ell_2+1).
\]
Therefore, all eigenstates of $ \Delta(H) $  form a basis of 
$ V^{(\ell_1)}_{+} \otimes V^{(\ell_2)}_{-}. $ In the present case 
we have thus proved the decomposition (\ref{CGdecomp}). 

%
%
\subsection{Odd and even dimensional representations}
\label{OE}

  We consider representations in the space $ V^{(\ell_1)} \otimes V^{(\ell_2)}_{\pm}, $ 
where the first (resp. second) space in the tensor product is 
odd (resp. even) dimensional. 
The CGC for this case can be computed in the same way as in the 
previous subsection. We here list some corresponding formulae and 
omit the computational detail. 
A highest weight state in $ V^{(\ell_1)} \otimes V^{(\ell_2)}_{\pm} $ 
has the form of (\ref{ell-Lambda}) where the summation variables run 
under the constraint $ m_1 + m_2 = \ell. $ 
Since $ m_1 $ (resp. $ m_2 $) is a integer (resp. an half-integer), 
$ \ell $ takes a half-integral value and the constant $ \eta $ remains 
in the expression of the weight:
\beq
  \Delta(H)\, e^{\ell}_{\ell}(\ell_1,\ell_2,\Lambda)
  = 
  \hf (\ell \pm \eta)\; e^{\ell}_{\ell}(\ell_1,\ell_2,\Lambda).
  \label{DHl-OE}
\eeq
The highest weight condition determines the coefficient $ C_{m_1,m_2}:$ 
\bea
  C_{m_1,m_2} 
  &=&  (-1)^{(\ell_1-m_1) \lambda + \hf (\ell_1-m_1) (\ell_1-m_1-1)}
  q^{\hf (\ell+1\pm \eta)(\ell_1-m_1)}
  \nn \\
  &\times& 
  \left(
    \frac{1}{\ku{2\ell_1}!}
    \frac{ \ku{\ell_1+\ell_2-\ell}! \ku{\ell_1+m_1}! \ku{\ell_2+m_2}! }
         { \ku{-\ell_1+\ell_2+\ell}! \ku{\ell_1-m_1}! \ku{\ell_2-m_2}!}
  \right)^{1/2} C_{\ell_1,\ell-\ell_1}.
  \label{C-sol-OE}
\eea
The factor $ C_{\ell_1,\ell-\ell_1} $ may be determined by normalization 
of the highest weight states. 
Other states in $ V^{(\ell_1)} \otimes V^{(\ell_2)}_{\pm} $  
are obtained by repeated applications of $ \Delta(V_-) $ on 
the highest weight states. The CGC for this case may be read off 
from the expression of state vectors:
\bea
  \CGC{\ell_1\ \ell_2\ \;\ell}{m_1\,m_2\,m} 
  &=& N_3(\ell_1,\ell_2,\ell,m ; q) 
   \sum_{m_1+m_2=m} (-1)^{(\ell_1-m_1) \lambda + \hf (\ell_1-m_1) (\ell_1-m_1-1)}   
  \nn \\
  &\times&
      q^{-\hf m_1(m+1 \pm \eta)}
    \left(  
      \frac{ \ku{\ell_1-m_1}! \ku{\ell_2-m_2}! }{ \ku{\ell_1+m_1}! \ku{\ell_2+m_2}! }
    \right)^{1/2}
  \nn \\
  &\times& \sum_k 
   (-1)^{k(\ell_1+\ell_2-m) + \hf k(k-1)} 
   q^{\hf k (\ell+m+1)}
  \nn \\
  &\times& 
   \frac{ \ku{\ell_1+\ell-m_2-k}! \ku{\ell_2+m_2+k}! }
        { \ku{\ell_1-\ell+m_2+k}! \ku{\ell_2-m_2-k}! \ku{\ell-m-k}! \ku{k}! }.
  \label{CGC_OE}
\eea
The factor that includes the sum over $k$ in (\ref{CGC_OE}) may be 
converted into identical form as (\ref{sumk_OO}). 
Thus the polynomial part of the CGC in (\ref{CGC_OE}) leads to 
the $Q$-Hahn polynomial with $ Q = -q $ and the parameters listed 
in (\ref{Hahn_label}). The same discussion as in the previous subsection  
completes the proof of the decomposition (\ref{CGdecomp}).

%
%
%
\setcounter{equation}{0}
\section{Even dimensional representations of \bm{$OSp_q(1/2)$} 
and little \bm{$Q$}-Jacobi polynomials}
\label{RepA}

  In this section, we compute even dimensional representations of 
the quantum supergroup $ {\cal A} \equiv OSp_q(1/2) $ 
by choosing a different basis set, and adopting a different method from 
\cite{Zou03}. We remark that precise theory of matrix representations 
of quantum group has been developed in \cite{Wor}, and 
that the odd dimensional representations of $ \A $ 
have been obtained in \cite{Zou03} and \cite{ACNS}. 
In \cite{Zou03}, an algebra $ {\cal A}(\sigma) $ generated by  
$ 2 \times 2 $ quantum supermatrix 
is set in the beginning and later its dual algebra is identified with $ \U. $ 
The representations of the algebra $ {\cal A}(\sigma)$ are obtained in a way 
parallel to \cite{MMNNU}. 
In contrast to this approach, we start with the algebra $\U,$ and then 
determine its dual basis. It follows the construction of the universal 
$ {\cal T}$-matrix and the representations of $\A$ are readily obtained 
by taking matrix elements of the universal ${\cal T}$-matrix in the 
representation space of the algebra $ \U. $ In addition to the easy and 
clear mechanism of our construction, the use of the universal 
${\cal T}$-matrix imparts following advantages: (i) Algebraic structure 
of  $ \A $ is made transparent in the construction of the universal 
${\cal T}$-matrix, as its basis set is determined explicitly. 
(ii) Nontrivial contribution of parity odd elements of $ \A $ to 
representations can be easily read off from the form of the universal 
${\cal T}$-matrix, that is, distinction from Lie supergroup $OSp(1/2)$ 
is emphasized. 

We divide this section in two parts. The first part contains a summary 
of the basis of the algebra $ \A$, and we also quote the universal 
${\cal T}$-matrix from \cite{ACNS}, where its detailed construction  
employing the Hopf duality between the ${\cal U}$ and ${\cal A}$ 
algebras is given. The second part is devoted to computation of the 
even dimensional representations of $ \A$ and their relation to the 
little $Q$-Jacobi polynomials.

%
%
\subsection{$ OSp_q(1/2)$ and Universal ${\cal T}$-matrix}
\label{OSPuT}

 The algebra $ \A, $ introduced in \cite{Kul,KR}, 
is a Hopf algebra dual to the algebra $\U.$ 
Two Hopf algebras $\U$ and $\A$ are in duality if there 
exists a doubly-nondegenerate bilinear form
$
\langle \; ,\; \rangle: \A \otimes \U 
\ \rightarrow \  {\mathbb C} 
$
such that, for $({\sf a}, {\sf b}) \in \A, 
({\sf u}, {\sf v}) \in \U$,
\bea
&&\langle {\sf a}, {\sf u v}\rangle = \langle \Delta_{\A}({\sf a}),
{\sf u} \otimes {\sf v}\rangle, \qquad
\langle {\sf a b}, {\sf u}\rangle = \langle {\sf a}\otimes {\sf b},
\Delta_{\U}({\sf u})\rangle,\nn\\
&&\langle {\sf a}, 1_{\U}\rangle = \epsilon_{\A}({\sf a}),\quad 
\langle 1_{\A}, {\sf u}\rangle = \epsilon_{\U}({\sf u}),\quad 
\langle {\sf a}, S_{\U}({\sf u})\rangle =  
\langle S_{\A}({\sf a}), {\sf u}\rangle.   
\label{dual_def}
\eea
The algebra $ \A $ is generated by  three elements, which are dual to 
the generators of the algebra $ \U$: 
\beq
  \langle x, V_+\rangle = 1, \quad
  \langle z, H\rangle = 1, \quad
  \langle y, V_-\rangle = 1. 
  \label{xyz_def}
\eeq
Thus, $x$ and $y$ are of odd parity, while $z$ is even. 
The generating elements satisfy the commutation relations:
\beq
  \{ x, y \} = 0, \qquad [z, x] = 2 \ln q\,\,x, 
  \qquad [z, y] = 2 \ln q\,\, y.
  \label{Acomm}
\eeq
Let the ordered monomials $E_{k \ell m} = V_{+}^{k} 
H^{\ell} V_{-}^{m},\,\, (k,\ell,m) \in (0, 1, 2,\cdots)$ be the 
basis elements of the algebra $\U.$ 
The basis elements $e^{k \ell m}$ of the dual Hopf algebra $\A$  
follow the relation
\beq
  \langle e^{k \ell m}, E_{k'\ell'm'} \rangle = 
  \delta^k_{k'} \delta^{\ell}_{\ell'} \delta^m_{m'}.
  \label{pair}
\eeq
The generating elements of the algebra $\A$ may be 
identified as $x = e^{1 0 0}, y = e^{0 0 1}$ and $z = e^{0 1 0}.$ 
The basis elements $e^{k \ell m}$ are ordered polynomials in the generating 
elements:
\beq
   e^{k \ell m} = \frac{x^k}{\ku{k}!} 
   \frac{(z+(k-m) \ln q)^{\ell}}{\ell !} \frac{ y^m}{\bra{m}!}.
   \label{enrs}
\eeq
Using the duality structure full Hopf structure of the algebra $ \A $  
has been obtained in \cite{ACNS}. We, however, do not list them here as 
it is not used in the subsequent discussions. 

 The notion of universal ${\cal T}$-matrix is a key feature capping 
the Hopf duality structure. Universal $ {\cal T}$-matrix 
for the superalgebra is defined by 
\beq
   {\cal T}_{e,E} = 
   \sum_{k\ell m} (-1)^{\widehat{e^{k\ell m}} (\widehat{e^{k\ell m}}-1)/2}\, 
   e^{k\ell m} \otimes E_{k\ell m},
   \label{uniTdef}
\eeq
where the parity of basis elements is same for two 
Hopf algebras $\U $ and $\A$:
\beq
   \widehat{e^{k\ell m}} = \widehat{E_{k\ell m}} = k + m.   \label{paritygen}
\eeq
Consequently, the duality relations 
(\ref{dual_def}) may be concisely expressed \cite{FG} in terms of the 
${\cal T}$-matrix as
\bea
    & & \uT{e,E} \uT{e',E} = \uT{\Delta(e),E}, \qquad 
        \uT{e,E} \uT{e,E'} = \uT{e,\Delta(E)}, \nn \\
    & & \uT{\epsilon(e),E} = \uT{e,\epsilon(E)} = 1, \quad \ \;
        \uT{S(e),E} = \uT{e,S(E)}.  \label{uniTprop}
\eea
where $e$ and $e^{\prime}$\, (resp. $E$ and $E^{\prime}$) refer to the 
two identical copies of algebra $\A$ (resp. $\U$). 

  Our explicit listing of the complete set of dual basis elements in 
(\ref{enrs}) allows us to obtain the universal ${\cal T}$-matrix as 
an operator valued function in a closed form:
\bea
  \uT{e,E} &=& 
  \left( \sum_{k=0}^{\infty} \frac{(x \otimes V_+q^H)^k}{\sq{k}!} \right) 
  \exp(z\otimes H) 
  \left(
     \sum_{m=0}^{\infty} \frac{(y \otimes q^{-H} V_-)^m}{\sqi{m}!} 
  \right)
  \nn \\
  &\equiv& {}_{\times}^{\times} 
  {\cal E}{\rm xp}_q (x \otimes V_+q^H) \; \exp(z \otimes H) \;
  {\cal E}{\rm xp}_{q^{-1}} (y \otimes q^{-H}V_-){}_{\times}^{\times},
  \label{uniTclosed}
\eea
where we have introduced a deformed exponential that is characteristic 
of the quantum $OSp_{q}(1/2)$ supergroup:
\beq
  {\cal E}{\rm xp}_q({\cal X}) \equiv \sum_{n=0}^{\infty} 
  \frac{{\cal X}^n}{\sq{n}!},
  \qquad
  \sq{n} = \frac{ 1 - (-1)^n q^n}{1 + q}.
  \label{deformedExp}
\eeq
The operator ordering has been explicitly indicated in (\ref{uniTclosed}). 
In \cite{DKLS}, using the Gauss
decomposition of the fundamental representation a universal 
$\cal T$-matrix for $ \U $ is given in terms of the 
standard $q$-exponential instead of the deformed exponential 
(\ref{deformedExp}) characterizing quantum supergroup. 
In the classical $ q \rightarrow 1 $ limit the universal ${\cal T}$-matrix 
(\ref{uniTclosed}) yields \cite{ACNS} the group element of the undeformed 
supergroup $ OSp(1/2)$. As the nilpotency relations $x^{2} = 0,\;
y^{2} = 0$ hold in the classical regime, we assume the following finite 
limits: 
\beq
\lim_{q \rightarrow 1} \frac{x^{2}}{q - 1} = \mathfrak{x},\qquad
\lim_{q \rightarrow 1} \frac{y^{2}}{q^{- 1} - 1} = \mathfrak{y}.
\label{xylim}
\eeq
It then  follows that in this limit the universal ${\cal T}$ matrix 
(\ref{uniTclosed}) reduces to an element of the classical supergroup 
$OSp(1/2)$: 
\beq
 {\cal G} = (1 \otimes 1 + x \otimes V_+ )\, 
 \exp (\mathfrak{x} \otimes V_{+}^{2}) \,\exp(z \otimes H) 
 \exp (\mathfrak{y} \otimes V_{-}^{2}) \, 
 (1\otimes 1 + y \otimes V_-).
 \label{Gclassical}
\eeq
The well-known existence of the classical $SL(2)$ subgroup structure 
generated by the elements $(V_{\pm}^{2}, H)$ of the undeformed 
$osp(1/2)$ algebra is evident from (\ref{Gclassical}). In fact, the 
correct limiting structure (\ref{Gclassical}) emphasizes that the 
quantum universal ${\cal T}$ matrix embodies the duality between the 
${\cal U}$ and the ${\cal A}$ algebras in a way that runs parallel to 
the familiar dual kinship between the classical Lie algebras and groups. 
Mappings from the universal ${\cal T}$-matrix to universal 
${\cal R}$-matrix of $\U$ also exist \cite{ACNS}.

%
%
\subsection{Representation Matrices}
\label{OSp_even}

  The closed form of the universal ${\cal T}$-matrix in 
(\ref{uniTclosed}) can be used to compute the representation matrices 
of the quantum supergroup $\A$ as has been done in \cite{ACNS}.  
To be explicit, we construct the representations of $\A$ by evaluating  
the matrix elements of the universal ${\cal T}$-matrix on $ V^{(\ell)}$ 
defined in \S \ref{UqRep}:
\bea
  \Tm{\ell}{m'm}(\lambda) &=& (\e{\ell}{m'}(\lambda),\uT{e,E}\, 
  \e{\ell}{m}(\lambda)) \nn \\
  &=& \sum_{abc} (-1)^{(a+c)(a+c-1)/2+(a+c)(\ell-m'+\lambda)} e^{abc}\,
  (\e{\ell}{m'}(\lambda), E_{abc}\, \e{\ell}{m}(\lambda)).
  \label{Tmatdef}
\eea
Assuming the completeness of the basis vectors 
$ \e{\ell}{m}(\lambda), $ we may verify the properties 
\beq
  \Delta(\Tm{\ell}{m'm}(\lambda)) = \sum_k \Tm{\ell}{m'k}(\lambda) \otimes 
  \Tm{\ell}{km}(\lambda),
  \qquad
  \epsilon(\Tm{\ell}{m'm}(\lambda)) = \delta_{m'm},
  \label{corepcond}
\eeq
which imply that the matrix elements (\ref{Tmatdef}) 
satisfy the axiom of comodule \cite{KS}. We may, therefore, regard 
$ \Tm{\ell}{m'm}(\lambda) $ as the 
$ 2\ell+1 $ dimensional matrix representation of the algebra $\A.$ 
We proceed to compute the matrix elements for the even dimensional 
representations (\ref{RepU_even}) in the present section. 
The explicit listing of the basis elements of $\A$ in (\ref{enrs}) 
renders the computation of the matrix elements straightforward. 
The computation is carried out by using two identities 
obtained by repeated use of (\ref{RepU_even}):
\beq
  V_+^a\, \ee{m} = (\pm 1)^a  
  \left(
    \frac{1}{ \ku{2}^a } \frac{ \ku{\ell-m}! }{ \ku{\ell+m}! }
    \frac{ \ku{\ell+m+a}! }{ \ku{\ell-m-a}! }
  \right)^{1/2}
  \ee{m+a},
  \label{Vpa}
\eeq
and
\beq
  V_-^c\, \ee{m} 
  = i^c (-1)^{c(\ell-m)+c(c-1)/2} 
  \left(
    \frac{1}{ \ku{2}^c } \frac{ \ku{\ell+m}! }{ \ku{\ell-m}! }
    \frac{ \ku{\ell-m+c}! }{ \ku{\ell+m-c}! }
  \right)^{1/2}
  \ee{m-c}. \label{Vmc} 
\eeq
We just quote the final result below:
\bea
  & & \Tm{\ell}{m'm} = (\pm 1)^{m'-m}
   (-1)^{\hf(m'-m)(m'-m-1)  + (m'-m)(\ell-m'+\lambda)} 
   e^{\pm  \frac{\pi}{4}(m'-m)i}
  \nn \\
  & & \quad \times \;
   q^{\hf m(m'-m)}
   \left(
      \frac{1}{ \ku{2}^{m'-m} }
      \frac{ \ku{\ell+m}! \ku{\ell+m'}! }{ \ku{\ell-m} ! \ku{\ell-m'}! }
   \right)^{1/2}
   \nn \\
   & & \quad \times \;
    \sum_c (\pm i)^{c}  (-1)^{c(\ell-m-1)} \,
    \frac{q^{-\hf c(m'-m)}}{ \ku{2}^c }
    \frac{ \ku{\ell-m+c}! }{ \ku{\ell+m-c}! }
    \nn \\
    & & \quad \times \;
    \frac{x^{m'-m+c}}{ \ku{m'-m+c}! }\, 
    \exp\left(
          \frac{m-c}{2}z \pm \frac{\pi i}{4 \ln q} z
        \right) \,
    \frac{y^c}{ \ku{c}! },
    \label{T_OSP_mat}
\eea
where the index $c$ runs over all non-negative integers maintaining the 
argument of $\ku{x}$ non-negative. 

  We now turn our attention to the polynomial structure built into the 
general matrix element (\ref{T_OSP_mat}) in terms of the variable 
\beq
  \zeta = -\frac{q^{-1/2}}{\ku{2}} xe^{-z/2}y.  \label{zetadef}
\eeq
We note that the variable (\ref{zetadef}) differs in sign from the 
corresponding one for the odd dimensional case. 
To demonstrate this, the product of generators in (\ref{T_OSP_mat}) for 
the case $m'-m \geq 0$ may be rearranged as follows:
\bea
  & & x^{m'-m+c} \, 
      \exp\left(
            \frac{m-c}{2}z \pm \frac{\pi i}{4 \ln q} z
          \right) \,
      y^c
    \nn \\
   & & \qquad 
     = (\mp i)^c q^{-mc} x^{m'-m} 
      \exp\left(
            \frac{m}{2}z \pm \frac{\pi i}{4 \ln q} z
          \right) \,
      x^c e^{-cz/2} y^c. \nn
\eea
The matrix element $\Tm{\ell}{m'm}(\lambda)$ may now be succinctly
expressed as a polynomial structure given below: 
\bea
  & & \Tm{\ell}{m'm} = (\pm 1)^{m'-m} (-1)^{(m'-m)(m'-m-1)/2 + (m'-m)(\ell-m'+\lambda)}
  \nn \\
  & & \quad \times \;
   e^{\pm  \frac{\pi}{4} (m'-m)i}\,
   \frac{q^{\hf m(m'-m)}}{ \ku{m'-m}!   }
   \left(
      \frac{1}{\ku{2}^{m'-m}}
      \frac{ \ku{\ell-m}! \ku{\ell+m'}! }{ \ku{\ell+m} ! \ku{\ell-m'}! }
   \right)^{1/2}
   \nn \\
    & & \quad \times \; 
    x^{m'-m} \,
    \exp\left(
          \frac{m}{2}z \pm \frac{\pi i}{4 \ln q} z
        \right) \,
    P^{\ell}_{m'm}(\zeta).
  \label{Tmatposi}
\eea
The polynomial $ \Pz{\ell}{m'm}$ in the variable $\zeta$ is defined by
\bea
  & & \Pz{\ell}{m'm} =
  \sum_c (-1)^{c(\ell-m)+ \hf c(c-1)} q^{-\hf c(m'+m-1)}
  \nn \\
  & & \hspace{2cm} \times
  \frac{\ku{m'-m}! \ku{l+m}! \ku{\ell-m+c}!}
   {\ku{m'-m+c}! \ku{\ell+m-c}! \ku{\ell-m}! \ku{c}!}\,
  \zeta^c,
  \label{Pzposi}
\eea
where the index $c$ runs over all non-negative integers maintaining the 
arguments of $\ku{x}$ non-negative. 
The polynomial (\ref{Pzposi}) is identical to the one appearing in 
the odd dimensional case \cite{ACNS}. 
For the case 
$ m'-m \leq 0,$ we make a replacement of the summation index $c$ with  
$a = m'-m+c$. Rearrangement of the generators now provides the 
following expression of the general matrix element:
\bea
  & & \Tm{\ell}{m'm} =  i^{m-m'}\,(-1)^{\hf(m-m')(m-m'+1)+(m-m')\lambda}
  \nn \\
  & & \quad \times \;
    e^{\mp \frac{\pi}{4}(m-m')i}
    \frac{q^{-\hf m'(m-m')}}{\ku{m-m'}!}
    \left(
       \frac{1}{\ku{2}^{m-m'}}
       \frac{ \ku{\ell+m} \ku{\ell-m'} }{ \ku{\ell-m} \ku{\ell+m'} }
    \right)^{1/2} 
   \nn \\
   & & \quad \times \; 
    \exp\left(
      \frac{m'}{2} z \pm \frac{\pi i}{4 \ln q} z
    \right) 
    y^{m-m'} P^{\ell}_{m'm}(\zeta),
  \label{Tmatnega}
\eea
where the polynomial $ \Pz{\ell}{m'm} $ for $ m'-m \leq 0 $ is defined by
\bea
  & & \Pz{\ell}{m'm} = 
  \sum_a (-1)^{a(\ell-m')+ \hf a(a-1)} q^{- \hf a(m'+m-1)}
  \nn \\
  & & \hspace{2cm} \times \frac{\ku{m-m'}! \ku{l+m'}! \ku{\ell-m'+a}!}
   {\ku{m-m'+a}! \ku{\ell+m'-a}! \ku{\ell-m'}! \ku{a}!}\,
  \zeta^a.
  \label{Pznega}
\eea
The polynomial (\ref{Pznega}) is also identical to the corresponding one 
obtained in \cite{ACNS}. 

  As seen above, the even dimensional representations of the 
algebra $\A$ have the same polynomial structure as the odd dimensional 
ones, though, between these two cases, the variable $ \zeta $ differs 
by a sign. Thus the polynomials (\ref{Pzposi}) and (\ref{Pznega}) are 
identified to the little $Q$-Jacobi polynomials with $ Q = -q.$ 
The little $Q$-Jacobi polynomials are defined via 
$ {}_2 \phi_1 $ \cite{GR}: 
\bea
  p^{(\alpha,\beta)}_m(z) 
  &=& {}_2 \phi_1(Q^{-m}, Q^{\alpha+\beta+m+1}; Q^{\alpha+1}; Q; Qz)
  \nn \\
  &=& 
  \sum_{n}  
  \frac{(Q^{-m};Q)_n (Q^{\alpha+\beta+m+1};Q)_n}{(Q^{\alpha+1};Q)_n (Q;Q)_n}
  (Qz)^n.
  \label{DefqJ2}
\eea
Rewriting our polynomials (\ref{Pzposi}) and (\ref{Pznega}) 
in terms of the shifted factorial with $ Q = -q, $ their 
identification is readily obtained.  
For the choice $ m'-m \geq 0,$ the polynomial structure reads 
\beq
  P^{\ell}_{m'm}(\zeta) =
     \sum_a  \,
     \frac{(Q^{-\ell-m};Q)_a\, (Q^{\ell-m+1};Q)_a}
          {(Q^{m'-m+1};Q)_a\, (Q;Q)_a} \,(Q \zeta)^a
     = p^{(m'-m,-m'-m)}_{\ell+m}(\zeta),
  \label{PtoJ1}
\eeq
and for the $ m'-m \leq 0$ case its identification is given by  
\beq
  \displaystyle
  P^{\ell}_{m'm}(\zeta) =
     \sum_a \,
     \frac{(Q^{-\ell-m'};Q)_a \, (Q^{\ell-m'+1};Q)_a}
          {(Q^{m-m'+1};Q)_a\, (Q;Q)_a} \,(Q \zeta)^a
     = p^{(m-m',-m'-m)}_{\ell+m'}(\zeta).
  \label{PtoJ2}
\eeq

  The even and odd dimensional representations of the algebra $\A$ 
have almost the same form. The fundamental difference of them is the 
factor 
$ 
  \exp( \pm \eta z/2)
$ 
appearing in the even dimensional representations. The factor is not 
well-defined in the classical limit of $ q \rightarrow 1.$ This feature of 
the even dimensional representations of the algebra $\A$ owes its genesis 
from the corresponding one of its dual algebra $\U.$ 

  To connect the results in this section to that of \cite{Zou03}, 
we consider the representation specified by $ \ell = \hf, \; \lambda = 0. $
We denote the matrix elements as follows:
\bea
  & & a = \Tm{\hf}{\hf \hf} 
        = \exp\left( \frac{z}{4} \pm \frac{\eta z}{2} \right)
          (1 + \zeta),
      \qquad
      b = \Tm{\hf}{\hf\; -\hf} = \pm \frac{ e^{\mp \frac{\pi}{4}i} }
                   { q^{\frac{1}{4}} \ku{2}^{\hf} } \, x d,
      \nn \\
   & & c = \Tm{\hf}{-\hf \hf} = i \frac{ e^{\mp \frac{\pi}{4}i} }
                   { q^{\frac{1}{4}} \ku{2}^{\hf} } \,d y,
       \qquad \qquad
       d = \Tm{\hf}{-\hf\; -\hf} = \exp\left( -\frac{z}{4} \pm \frac{\eta z}{2}  \right).
       \label{T_half}
\eea

The commutation relations satisfied by the matrix elements may be
immediately derived: 
\bea
  & & ab = \pm i q^{\hf} ba, \qquad
      ac = \pm i q^{\hf} ca, \qquad
      bc = - cb,
  \nn \\
  & & bd = \mp i q^{\hf} db, \qquad 
      cd = \mp i q^{\hf} dc, \qquad
      [a,\, d] = -(1+q) bc.
  \label{comm_half}
\eea

The central element $ ad + q bc $ commutes with $ a, d $ and 
anti-commutes with $ b, c.$ 
Thus the representations specified by $ \ell = \hf,\; \lambda = 0 $ 
is precisely same as the algebra $ {\cal A}(\sigma) $ used in 
\cite{Zou03}.

%
%
%
\section{Even dimensional Covariant Spaces of $\bm{OSp_q(1/2)}$}
\label{CovSpa}
\setcounter{equation}{0}

  It has been observed \cite{AC,AC2} that noncommutative spaces 
covariant under the action of a finite dimensional representation of 
quantum groups such as $ SL_q(2)$ or $ OSp_q(1/2)$ may be obtained by 
using the CGC. The method developed in \cite{AC,AC2} is outlined below. 
We introduce an algebraic structure on a given representation space 
$ V^{(\ell)}. $ Namely, assuming a multiplication map 
$ \mu : V^{(\ell)} \otimes V^{(\ell)} \ \rightarrow \ V^{(\ell)}, $ 
we determine a consistent set of commutation relations among bases of 
$ V^{(\ell)}$ that may be regarded as the generators of noncommutative 
spaces. Specifically, for the highest weight $ \ell $ representation 
of $ OSp_q(1/2), $ we construct the following composite object:
\beq
 E^{L}_M(\Lambda) =
     \sum_{m_1,m_2\atop m_{1}+m_{2} = M} 
    \CGC{\ell\quad \ell\ \;\,L}{m_1\,m_2\,M}\,\,
    \e{\ell}{m_1}(\lambda)  \e{\ell}{m_2}(\lambda),
  \label{compobj}
\eeq
where $ \Lambda = 2\, \ell - L $ (mod 2), where $\ell$ is of integral or 
half-integral value. Then it may be proved that the following relations 
are covariant under the right coaction of the highest weight $ \ell $ 
representation of $ OSp_q(1/2): $
\bea
  & & E^0_0(0) = r, \label{E00} \\
  & & E^{\ell}_M(\Lambda) = \xi e^{\ell}_M(\lambda), \label{Eell} \\
  & & E^L_M(\Lambda) = 0, \quad (L \neq \ell, 0), \label{EL}
\eea
where $r$ and $ \xi $ are parameters. In the $ q \rightarrow 1 $ limit 
$ \xi \rightarrow 0, $ and $ \xi $ is regarded as a Grassmann number
if the parity of the two sets of vectors in (\ref{Eell}) differ: 
$ \Lambda \neq \lambda $ (mod 2).

Although we have obtained a set of covariant commutation properties, 
the simultaneous use of all relations from (\ref{E00}) to (\ref{EL}) 
gives an inconsistent result, since some of them do not have correct 
classical limits. In order to obtain a consistent covariant algebra, we 
have to make a choice regarding the relations to be used for defining 
the algebra. Then their consistency has to be verified. The consistency 
requirements are listed below:
\begin{enumerate}
\renewcommand{\labelenumi}{(\alph{enumi})}
 \item The constant $r$ commutes with all generators.
 \item The associativity of products of generators need to be 
       maintained. 
\end{enumerate}

  Employing the above procedure in conjunction with the CGC given in 
(\ref{EE}), we now construct covariant noncommutative spaces of 
dimensions two and four. 

\medskip\noindent
Case 1: $ \ell = \hf $ 

 The allowed values of $ L $ are $ 0 $ and $ 1. $ 
We rewrite the basis of $ V^{(1/2)} $ as follows:
\[
 \e{1/2}{1/2}(\lambda) \ \rightarrow \ x, \qquad
 \e{1/2}{-1/2}(\lambda) \ \rightarrow \ y.
\]
We first consider the case of  $ \lambda = 0, $ when 
$ x\;\; (\hbox{resp.}\; y) $ is of 
even (resp. odd) parity. Covariant relations for $ L = 1 $ obtained from 
(\ref{EL}) read
\beq
  x^2 = y^2 = 0, \qquad q^{1/4} xy + q^{-1/4} yx = 0. \label{L1-1} 
\eeq
These relations are unacceptable as a definition of the covariant 
noncommutative space as there the even element $x$ becomes nilpotent. We 
rather regard $ L = 0 $ relation obtained from (\ref{E00}) as a 
definition of the covariant space:
\beq
  xy + q^{1/2} yx = r. \label{L0-1}
\eeq
We now illustrate $\lambda = 1$ case where $x\;\;(\hbox{resp.}\; y)$ 
is of odd (resp. even) parity. One can see that $ L = 1 $ relations are 
rejected again by the same reason. We thus obtain a covariant space 
from (\ref{E00}):
\beq
  xy - q^{1/2} yx = r. \label{L0-2}
\eeq
Setting $ r = 0 $ in (\ref{L0-1}) or (\ref{L0-2}), the quantum superspaces 
found in literatures (\textit{e.g.} \cite{Ma,KU}) are recovered. However, 
(\ref{L0-1}) or (\ref{L0-2}) gives the most general two-dimensional 
covariant superspaces.

\medskip\noindent
Case 2: $ \ell = \frac{3}{2} $ 

  The index $ L $ ranges the integral values from 3 to 0. 
We rewrite the basis of $ V^{(3/2)} $ as
\[
  \e{3/2}{3/2}(\lambda) \ \rightarrow \ x, \quad
  \e{3/2}{1/2}(\lambda) \ \rightarrow \ y, \quad
  \e{3/2}{-1/2}(\lambda) \ \rightarrow \ z, \quad
  \e{3/2}{-3/2}(\lambda) \ \rightarrow \ w.
\]
We study the case of $ \lambda = 0, $ since the example $ \lambda = 1 $ 
yields almost identical results except for some sign differences. For 
the choice $ \lambda = 0, $ the generating elements 
$ x, z \;\;(\hbox{resp.}\; y, w) $ are of even (resp. odd) parity. The 
results corresponding to $ L = 3 $ obtained from (\ref{EL}) contain an 
unacceptable relation: $ x^2 = 0. $ We, therefore, discard these 
equations. The equations for the $ L = 1, 2 $ cases obtained from 
(\ref{EL}) provide six commutation relations and two additional constraints:
\bea
  & & xy + q^{3/2} yx = 0, \qquad 
      xz - q^3 zx = 0,
  \nn \\
  & & (q^{-2}-1+q) xw - \ku{2} wx + q^{-1/2}(q^2-1+q^{-2}) yz = 0,
  \nn \\
  & & (q^2-1+q^{-1}) yz + q \ku{2} zy + q^{1/2} \ku{3} xw = 0, 
  \label{4Drel1} \\
  & & yw - q^3 wy = 0, \qquad
      zw + q^{3/2} wz = 0, \nn
\eea
and
\beq
  y^2 = q^{-3/2} \sqrt{\ku{3}} \, xz, \qquad
  z^2 = q^{-3/2} \sqrt{\ku{3}} \, yw. 
  \label{4Drel2}
\eeq
It turns out that these relations do not satisfy the consistency 
condition (b). For instance, two ways of reversing $ xyz $ to $ zyx $ 
do not give identical result. We thus incorporate the relation 
corresponding to $ L = 0$. We regard this relation as an additional 
constraint after setting $ r = 0 $ in (\ref{E00}), and thereby make the 
four dimensional covariant space well-defined. The constraint reads 
$ yz = -q^{-3/2} \ku{3} xw. $ Employing this constraint we may simplify 
the commutation properties (\ref{4Drel1}) and (\ref{4Drel2}). The  
requirement that the simplified relations obey the consistency condition 
(b) is also verified. We, therefore, introduce the four dimensional 
covariant space defined by the six commutation relations 
\bea
  & & x y = - q^{3/2} yx, \qquad 
      x z = q^3 z x, \qquad xw = -q^{9/2} wx,
  \nn \\
  & & y z = -q^{3/2} z y, \qquad 
      y w = q^3 w y, \qquad zw = -q^{3/2} wz,
  \label{4Ddef1}
\eea
and three constraints
\beq
  y^2 = q^{-3/2} \sqrt{\ku{3}} \, xz, \qquad
  z^2 = q^{-3/2} \sqrt{\ku{3}} \, yw, \qquad
  yz = -q^{-3/2} \ku{3} xw.
  \label{4Ddef2}
\eeq

%
%
%
\section{Concluding Remarks}
\label{ConclR}

  We have seen intimate relations between the representations of 
the algebras $ \U, \A $ and basic hypergeometric functions. 
Existence of even dimensional representations makes the representation 
theory of the quantized $osp(1/2)$  algebra richer and more interesting 
than the one of the classical Lie superalgebra $ osp(1/2).$ Especially, 
the fundamental representation of $ osp(1/2)$ is three dimensional, 
while the corresponding representation of $ \U $ can be further 
decomposed into product of two dimensional ones. In other words, 
quasi-particles described by the three dimensional representation 
of $ \U $ can be regarded as a composite of more fundamental objects. 
Such situation, hopefully, may be realized in some physical models. 
A byproduct of the even dimensional representations of the 
$osp_{q}(1/2)$ algebra is that new noncommutative spaces covariant under 
the coaction of the quantum group $OSp_{q}(1/2)$ may be constructed 
via the Clebsch-Gordan decomposition. The representations of these 
noncommutative spaces for root of unity values of $q$, for instance, may 
be relevant for some physical problems.

  Turning to the representations of Lie superalgebras, 
little seems to be known about their relations to hypergeometric 
functions. This may be explained by the appearance of the $ Q = -q $ 
polynomials for the case of the algebras $ \U $ and $ \A. $ 
The classical limit of such polynomials has somewhat complicated structure 
as they have to be evaluated at $ Q = -1.$  It may be difficult to find 
such polynomials starting from the representations of classical objects 
such as $ osp(1/2) $ and $ OSp(1/2).$  In this sense, study of the 
representations of the quantum superalgebras gives deeper understanding 
of the representation theory of the Lie superalgebras. 
It is known that there is a one-to-one correspondence between the finite 
dimensional representations of $ osp(1/2n) $ and $ so(2n+1)$ except for 
the spinorial ones. For quantum algebras, this is explained 
\cite{Zhang} by the isomorphism between $ U_q[osp(1/2n)] $ and 
$ U_{-q}[so(2n+1)]$, which holds on the non-spinorial representation 
spaces. Our work confirms that for the even dimensional representations 
for which the said isomorphism is not known are still characterized by 
the $Q = -q$ polynomials. This may be a more general feature of the 
quantum supergroups.

%
%
%
\section*{Acknowledgements}
We thank S. Blumen for introducing us to the work of Zou \cite{Zou}. 
We also thank the referee for useful comments. The work of N.A. 
is partially supported by the grants-in-aid from JSPS, Japan (Contract 
No. 15540132). Another author (S.S.N.M.) is supported by the University 
Grants Commission, Government of India. 

%
%
%

\end{document}